\documentclass[10pt,english,a4paper]{article}
\usepackage{a4wide}

\usepackage{amsmath}
\usepackage{amssymb}

\title{Conjectures on the ring of commuting matrices}
\author{Freyja Hreinsdottir}
\date{}

\begin{document}
\maketitle

\begin{abstract}
Let $X=(x_{ij})$ and $Y=(y_{ij})$ be generic $n$ by $n$ matrices and $Z=XY-YX$. Let
$S=k[x_{11},\ldots,x_{nn},y_{11},\ldots,y_{nn}]$, where $k$ is a field, let
$I$ be the ideal generated by the entries of $Z$ and let $R=S/I$. We give a conjecture on the first syzygies of $I$, show how these can be used to give a conjecture on the canonical module of $R$. Using this and the Hilbert series of $I$ we give a conjecture on the
Betti numbers of $I$ in the $4 \times 4$ case. We also give some guesses on the structure of the resolution in general.
\end{abstract}

\section{Introduction} \label{intro}
Throughout this article we let $R$ be the ring defined in the abstract.

\smallskip
\noindent
It was shown by Motzkin and Taussky~\cite{mot} that
the variety of commuting matrices in $M_n(k)$ is irreducible of dimension
$n^2+n$. Gerstenhaber~\cite{gerst} also showed that the variety is
irreducible. From this it follows that
$\operatorname{Rad}(I)$
is prime and that the dimension of $R$ is $n^2+n$.

\smallskip
\noindent
It has been conjectured that $R$ is Cohen-Macaulay and
this has been
shown for $n=3$ in~\cite{bss} and for $n=4$ in~\cite{frh1}. It has also
been conjectured that $R$ is a domain which follows from the ring
being CM (see~\cite{vasc}).

\smallskip
\noindent
Recently, Knutson~\cite{kn} proved that the off-diagonal elements in $XY-YX$
form a regular sequence.

\smallskip
\noindent
For the cases $n=2,3,4$ the computer programs {\sc Macaulay}~\cite{bast} and {\sc Macaulay 2}~\cite{m2} can be used to
compute a
Gr\"{o}bner basis (and thus the Hilbert series) of the ideal $I$. The resolution can be computed for
$n=2,3$ and partially for several $n\geq 4$. For the ideal generated by the off-diagonal 
elements, which we call $J$, a Gr\"{o}bner basis can easily be computed for the cases $n=2,3$ but we have not yet found a term 
order that works for $n=4$. In this article we use {\sc Macaulay} and {\sc Macaulay 2} to verify 
many conjectures. 

\smallskip
\noindent
By exploring simple facts concerning the trace of a matrix we can give "many" first syzygies of the
ideal $I$ and in section $2$ we give a conjecture on the first Betti numbers. In sections $3$
and $4$ we use our syzygy conjecture to give a conjecture on the generators of the ideal $(J:I)$ and
the canonical module of $S/I$. We then
compute a partial resolution of the (conjectured) canonical module in the case $n=4$ and splice this together with a partial resolution of $I$ to give a conjecture of the Betti numbers of $I$ in that case.

\smallskip
\noindent
In section $5$ we give some guesses (mostly based on computer calculations) on the resolution in general
and in section $6$ we comment on Knutson's conjecture concerning the prime ideals of $J$.

\section{First syzygies} \label{syzy}

We restate here a conjecture on the first syzygies that was given in~\cite{frh4}.

Write $I=(f_1, \ldots, f_{n^2})$, with $f_1=Z_{11}$,
$f_2=Z_{21}$, \ldots, $f_{n^2}=Z_{nn}$, where $Z=XY-YX$. A syzygy on $I$ is an $n^2$-tuple
$(a_1,\ldots,a_{n^2})$ such that
\begin{equation} \label{releq2}
f_1a_1 + f_2a_2 + \cdots +f_{n^2}a_{n^2} =0.
\end{equation}
This can be rewritten as
\begin{equation} \label{treq1}
\operatorname{tr} \left(
\begin{bmatrix}
a_{1}&\cdots&a_{n}\\
a_{n+1} & \cdots& \vdots  \\
\vdots& &\vdots\\
a_{n^2+n-1}&\ldots&a_{n^2}
\end{bmatrix}
\begin{bmatrix}
f_{1}&\cdots&f_{n^2-n+1}\\
f_2  & \cdots& \vdots  \\
\vdots& &\vdots\\
f_n&\ldots&f_{n^2}
\end{bmatrix}
\right)
=0
\end{equation}
i.e. as
\begin{equation} \label{treq2}
\operatorname{tr}(A(XY-YX))=0.
\end{equation}
So solving (\ref{treq2}) for $A$ is equivalent to solving (\ref{releq2}) for
$(a_1, \ldots, a_{n^2})$. We claim the following:

\medskip
\noindent
{\bf Conjecture:} the module of first syzygies on $I$ is generated by the Koszul relations of $I$
and matrices $A$ that are polynomials in $X$ and $Y$. The highest degree of a first syzygy is $n-1$ and
in degree $h$ we get generators of bidegrees $(h,0)$, $(h-1,1)$, \ldots, $(0,h)$ (considering the bidegree ($x$-deg, $y$-deg)).

\smallskip
\noindent
We do not have a proof of this conjecture but we show below that a number of solutions to (\ref{treq2}) exist and compare this
with computer calculations.

\smallskip
\noindent
In general we have $\operatorname{tr}(BC)=\operatorname{tr}(CB)$ for matrices $B$ and $C$
so we get that $\operatorname{tr}(A(XY-YX))=0$ whenever the matrix $A$ commutes with
$X$ or $Y$. This gives that any polynomial in $X$ and any polynomial
in $Y$ is a solution to equation~(\ref{treq2}).
We note the following:
\begin{quote}
If $M$ is a monomial in $X$ and $Y$ then $\operatorname{tr}(M(XY-YX))=0$
if $MXY$ can be cyclically permuted into $MYX$.
\end{quote}

\noindent
and

\begin{quote}
Let $M_1$ and $M_2$ be monomials in $X$ and $Y$ such that
$\operatorname{tr}(M_1(XY-YX)) \neq 0$ and $\operatorname{tr}(M_2(XY-YX)) \neq 0$. If $B=M_1+M_2$, then
$\operatorname{tr}(B(XY-YX))=0$
if $M_1XY$ can be cyclically permuted into $M_2YX$ and $M_2XY$ can
be cyclically permuted into $M_1YX$. If $B=M_1-M_2$ then $\operatorname{tr}(B(XY-YX))=0$ if
$M_1XY$ can be cyclically permuted into $M_2XY$ and $M_1YX$ can be cyclically
permuted into $M_2YX$
\end{quote}

\noindent
Using the above we can guess a number of solutions:
\begin{description}

\item[Degree $0$:]  Here we only have one syzygy $A=E$ (the identity matrix), i.e. the ideal
is minimally generated by $n^2-1$ elements.

\item[Degree $1$:] We have $\operatorname{tr}(X(XY-YX)) = \operatorname{tr}(X^2Y)-\operatorname{tr}(XYX) =0$ so
$A=X$ is a solution and similarly we get that $A=Y$ is a solution.
The two syzygies we get are obviously independent over $k$ as they have the bidegrees $(1,0)$ and $(0,1)$.
In~\cite{frh2} we proved that these are the only ones of degree $1$. 

\item[Degree $2$:] We see that $A=X^2$ and $A=Y^2$ are solutions.
The only other monomials in $X$ and $Y$ are $XY$ and $YX$ and
neither of those is a solution.  We have
\begin{eqnarray*}
\lefteqn{\operatorname{tr}((XY+YX)(XY-YX))} \\
 &=& \operatorname{tr}(XYXY)-\operatorname{tr}(XYYX)+
       \operatorname{tr}(YXXY)-\operatorname{tr}(YXYX)  \\
 &=& \operatorname{tr}(XYXY)-\operatorname{tr}(X^2Y^2)+\operatorname{tr}(X^2Y^2)-\operatorname{tr}(XYXY)           \\
              &=& 0
\end{eqnarray*}
so $A=XY+YX$ gives a syzygy. We thus have syzygies of bidegrees $(2,0),(1,1),(0,2)$.

\item[Degree $3$:] Here we get at least the monomial solutions $X^3$, $Y^3$, 
$XYX$, $YXY$ and the binomial solutions $X^2Y+YX^2$, $XY^2+Y^2X$. Macaulay calculations
indicate that it is enough to take one syzygy of each bidegree i.e. $X^3$, $Y^3$,$XYX$, $YXY$ will do.

\item[Degree $4$:] $X^4$, $Y^4$, $X^3Y+YX^3$, $Y^3X+XY^3$,$X^2YX+XYX^2$,  $Y^2XY+YXY^2$ and $XY^2X-YX^2Y$.

\item[Degree $5$:] $X^5$, $Y^5$, $X^2YX^2$, $Y^2XY^2$, $X^4Y+YX^4$, $XY^4+Y^4X$, $XYX^2Y+YX^2YX$, $YXY^2X+XY^2XY$.

\end{description}

\noindent
The syzygies given above work for any $n$. 
For $n=2,3,4$ (and partially for $n=5,6,7$) we can compare this with Macaulay calculations.

\subsection{$n=2$}

\begin{verbatim}

% betti s2
; total:   3     2
; -----------------
;     2:   3     2

\end{verbatim}
Which means that we have $3$ generators and the only syzygies we get are the $2$ linear ones.  The matrices $X$ and $Y$
are $2 \times 2$ matrices so they satisfy a characteristic polynomial
of degree $2$, i.e.
$X^2-\operatorname{tr}(X)X+\operatorname{det}(X)E=0$ so the syzygy given by $X^2$ can be written
in terms of smaller degree syzygies. Similarly for $Y^2$.
For $2 \times 2$
matrices (see e.g.~\cite{for}) we have the following identity:
\begin{equation} \label{traceid}
YX=(\operatorname{tr}(XY)-\operatorname{tr}(X)\operatorname{tr}(Y))E+
\operatorname{tr}(Y)X+\operatorname{tr}(X)Y-XY
\end{equation}
So the syzygy that $XY+YX$ gives can be written in terms of lower
degree syzygies.

\subsection{$n=3$}

We get the following Betti numbers

\begin{verbatim}

% betti s3
total:   8    33
-----------------
    2:   8     2
    3:   -    31

\end{verbatim}

\noindent
As expected we get $2$ linear first syzygies.
There are $31$ first syzygies of degree $2$, $\binom{8}{2}=28$ of those
are the trivial syzygies (Koszul relations) and the $3$ nontrivial ones correspond
to $A=X^2$, $A=Y^2$ and $A=XY+YX$.
There are no syzygies of degree $3$ so the solutions from before given by $X^3$, $Y^3$, $XYX$, $YXY$,
much be linear combinations of the syzygies of smaller degree. The characteristic equation
takes care of $X^3$ and $Y^3$ and we get

\begin{eqnarray*}
XYX  &=& \frac{1}{2}(x_1+x_5)(XY+YX) +y_9X^2  + (x_2x_4-x_1x_5)Y  \\
        & &  +(x_3y_7-x_1y_9+x_6y_8-x_5y_9)X +cE + \sum aT_a
\end{eqnarray*}
where the $T_a$ are matrices we get from the trivial syzygies.
The coefficient of $E$ is
\[
c=\operatorname{det}
\begin{pmatrix}
x_1 & x_2 & x_3 \\
x_4 & x_5 & x_6 \\
y_7 & y_8 & y_9
\end{pmatrix}
\]
and the coefficients of $XY+YX$, $X^2$, $X$ and $Y$ are given by traces and determinants
of minors of this matrix.

\noindent
Considering $n=4$ and $n=5,6,7$ (partial computation) we give the conjecture below on
the first Betti numbers. We use the notation of Macaulay 2 to display the Betti numbers, i.e. the number in column $i$ row $j$ (starting with column $0$, row $0$) 
is $\beta_{i,i+j}$.

\[
\begin{array}{rllcrr}
total:  &1 &  n^2-1  & \binom{n^2-1}{2}+\binom{n+1}{2} -1 \\
0: & 1 & .  &.& \\
1:  &.&   n^2-1  & 2 \\
2:  &.&    .    & \binom{n^2-1}{2} + 3  \\
3:  &.&    .    & 4     \\
4:  &.&    .    & 5     \\
.   &.&    .   & .     \\
.   &.&    .   & .     \\
.   &.&    .   & .     \\
n-1:  &.&     .   & n   \\
n:   &.&  .& . 
\end{array}
\]
Where the $\binom{n^2-1}{2}$ syzygies of degree $2$ are the Koszul relations.

We also conjecture that the $k$ syzygies of degree $k-1$ have the following
bidegrees (i.e. ($x$-deg, $y$-deg))
\[
(k-1,0), (k-2,1), \ldots, (1,k-2), (0,k-1).
\]

\section{The ideal quotient $(J:I)$}

In the following let $J$ be  the ideal generated by the off-diagonal elements of $XY-YX$.
Knutson~\cite{kn} has shown that these elements form a regular sequence for any $n$.
It is known~\cite{mot} that the
height of $I$ is $n^2-n$ which is equal to the number of generators of $J$
so these form a maximal regular sequence in $I$.

In this section we study the ideal $(J:I)$ and use our conjecture on the first syzygies
to give a conjecture on its generators. For the cases $n=2,3$ a Gr\"{o}bner basis of $J$ can be computed
using Macaulay so we can test the conjecture. We demonstrate first the cases $n=3,4$.

\subsection{$n=3$}
The nontrivial syzygies on $I$ are given by $A \in \{E,X,Y,X^2,Y^2,XY+YX\}$.

The ideal $I$ is generated by
$(f_1, \ldots ,f_9)$ where $f_1$, $f_5$ and $f_9$ are from the diagonal of
$XY-YX$ and $J=(f_2,f_3,f_4,f_6,f_7,f_8)$. Pick $3$ different syzygies, $A$, $B$ and $C$. Then

\begin{eqnarray*}
a_1f_1 +a_5f_5+a_9f_9 &=& a_2f_2+ a_3f_3 + a_4f_4+a_6f_6+a_7f_7+a_8f_8 \\
b_1f_1 +b_5f_5+b_9f_9 &=& b_2f_2+b_3f_3+b_4f_4+ b_6f_6 + b_7f_7+b_8f_8 \\
c_1f_1 +c_5f_5+c_9f_9 &=& c_2f_2+c_3f_3+c_4f_4+ c_6f_6 + c_7f_7+c_8f_8 \\
\end{eqnarray*}
so 
\[
\det \begin{bmatrix}
a_1 & b_1 & c_1\\
a_5 & b_5 & c_5\\
a_9 & b_9 & c_9
\end{bmatrix} \cdot f_i \in J \mbox{\quad for i=1, 5, 9}.
\] 
Direct calculations using {\sc Macaulay}
give that it suffices to take the generators of $J$ and the elements given by
$(A,B,C) \in \{ (E,X,Y),(E,X,X^2),(E,X,Y^2),(E,Y,Y^2),(E,Y,X^2) \}$
to get {\em all} the generators of $(J:I)$. The bidegrees of these additional generators
are $(1,1)$, $(3,0)$, $(2,1)$, $(1,2)$ and $(0,3)$. 

\noindent
We can also partially see directly that these suffice, denote by $A_d$ the diagonal of the matrix $A$ and by $1_d$ the
diagonal of $E$. Consider the $3 \times 4$ matrix $[1_d,A_d,B_d,C_d]$ and add one row by repeating say the first
row. We then have a $4 \times 4$ matrix whose determinant is zero and expanding by the first row we get
\[
0= \det[A_d,B_d,C_d] - A_{11}\det[1_d,B_d,C_d] + B_{11}\det[1_d,A_d,C_d]-C_{11}\det[1_d,A_d,B_d]
\]
So we only have to consider triples of the form $(E,A,B)$. There are $10$ of these, the ones given above and
$(E,X,XY+YX)$, $(E,Y,XY+YX)$, $(E,X^2,Y^2)$, $(E,X^2,XY+YX)$ and $(E,Y^2,XY+YX)$. We have
\[
\det[1_d,X_d,XY+YX_d] - 2\det[1_d,Y_d,X^2_d] - 2{\rm tr} X \det[1_d,X_d,Y_d] = x_2f_2-x_3f_3-x_4f_4+x_6f_6+x_7f_7-x_8f_8
\]
so the triple $(E,X,XY+YX)$ (and for symmetry reasons $(E,Y,XY+YX)$) does not give a new generator. It is probably possible to get similar simple equations explaining why $\det[1_d,X^2_d,Y^2_d]$, $\det[1_d,X^2_d,XY+YX_d]$ and
$\det[1_d,Y^2_d,XY+YX_d]$ are not needed as minimal generators of $(J:I)$.

\smallskip
\noindent
It seems that it suffices to use enough syzygies to
get one generator of each bidegree.

\subsection{$n=4$}
The nontrivial syzygies on $I$ are given by $A=X,Y,X^2,Y^2,XY+YX, X^3, XYX, YXY, Y^3$.

\smallskip
\noindent
Similarly to the case $n=3$ we pick $E$ and $3$ more syzygies and get that the determinant of a matrix consisting of the diagonals gives an element in $(J:I)$. Assuming that we get only one element in $(J:I)$ of each bidegree we get the following possibilites (ordered by total degree):

\bigskip
\begin{tabular}{lll}
total degree & triple of syzygies & bidegree of element in $(J:I)$ \\ \hline
$4$ & $(X,Y,X^2)$, $(X,Y,XY+YX)$, $(X,Y,Y^2)$  & $(3,1)$, $(2,2)$, $(1,3)$ \\
$5$ & $(X,Y,X^3)$, $(X,Y,XYX)$, $(X,Y,YXY)$, $(X,Y,Y^3)$ & $(4,1)$, $(3,2)$, $(2,3)$, $(1,4)$ \\
$6$ & $(X,X^2,X^3)$,$(Y,X^2,X^3)$, $(X,Y^2,X^3)$, $(X,X^2,Y^3)$   & $(6,0)$, $(5,1)$, $(4,2)$, $(3,3)$ \\
    & $(Y,X^2,Y^3)$, $(X,Y^2,Y^3)$, $(Y,Y^2,Y^3)$                 & $(2,4)$, $(1,5)$, $(0,6)$   \\
$7$ & $(XY+YX,X^2,X^3)$, \ldots                                    & $(6,1)$, \ldots  
\end{tabular}

\bigskip
\noindent
We believe that it is enough to take the elements in total degrees $4,5$ and $6$. So the highest degree
is high enough for us to pick syzygies in the $x$-variables only (and $y$-variables only).
This is partially based on 
comparison of the Hilbert series with a conjectured canonical module (see section $4$).

\subsection{General case}

We generalize the idea above for any $n$ and conjecture that $(J:I)$ is generated by
the elements of $J$ and elements of the form $u=\det U$ where $U$ is an $n \times n$ 
matrix whose colums are the diagonals of $E$ and the matrices defining the syzygies.
Below we give a conjecture on the degrees of these elements.

\smallskip
\noindent
We consider the following table of possible bidegrees of syzygies and find the smallest total degree of $u$
we can get from picking $n$ syzygies:

\medskip
\begin{center}
\begin{tabular}{ccccccc}
$(0,0)$ \\
$(1,0)$, $(0,1)$ \\
$(2,0)$, $(1,1)$, $(0,2)$ \\
$(3,0)$, $(2,1)$, $(1,2)$, $(0,3)$ \\
\ldots \ldots
\end{tabular}
\end{center}

\smallskip
\noindent
We pick bidegrees from the $k$ first rows where $k=\lfloor -\frac{1}{2} + \sqrt{2n+\frac{1}{4}} \rfloor$.
If $s:=n-\frac{k(k+1)}{2} \neq 0$  we pick $s$ bidegrees from row $k+1$ starting from the left. The total
bidegree from the first $k$ rows is $(\frac{k}{6}(k^2-1), \frac{k}{6}(k^2-1))=:(a,a)$. We have $2$ cases:

\medskip
\noindent
{\bf Case $s=0$:} then $n=\frac{k(k+1)}{2}$ and we get exactly $n$ matrices from the first $k$ rows. The smallest
possible total degree of $u=\det U$ is $d_{min}= 2\frac{k}{6}(k^2-1) = \frac{k}{3}(k^2-1)$ where
$k= -\frac{1}{2} + \sqrt{2n+\frac{1}{4}}$. As there is one possibility of picking the $n$ matrices we
get one generator of this minimal degree (the values of $n$ where this occurs are for instance $n=3, 6, 10, \ldots$).

\smallskip
\noindent
{\bf Case $s \neq 0$:} we pick $s$ bidegrees from row $k+1$ starting from the left. The total bidegree
is $(sk-\frac{s(s-1)}{2},\frac{s(s-1)}{2}):=(c-b,b)$. The smallest total degree is then:
\[
d_{min} = \frac{k}{3}(k^2-1) + sk
\] 
Considering the different possibilities of bidegrees of elements of this total degree (picking bidegrees further
to the right in the table of bidegrees) we get elements of bidegrees:
\[
(a+c-b,a+b), (a+c-b-1,a+b+1), \ldots, (a+b,a+c-b)
\]
the total number of such elements is $c-2b+1=s(k-s+1)+1$.

\medskip
\noindent
To get the number of elements in the next smallest degree we count the possibilities of picking $n$ bidegrees
and skipping $(1,0)$ or $(0,1)$ etc. In each total degree except the largest we get elements of bidegrees 
of the following form:
\[
(r,t), (r-1,t+1), \ldots, (r-(r-t),t+(r-t))
\]
where $r \neq 0$ and $t \neq 0$.

\medskip
\noindent
We believe that the largest total degree we need for a generator of $(J:I)$ is the one where we pick $n$ matrices of bidgrees
\[
(0,0), (1,0), \ldots, (n-1,0)
\]
giving us the total degree
\[
d_{max} = 0+ 1 + \cdots + (n-1) = \frac{n(n-1)}{2}.
\]
We get $\frac{n(n-1)}{2} +1 $ elements of this total degree, one for each possible bidegree.

\medskip
\noindent
In section $5$ we will see that this agrees with some guesses we have on the Betti numbers of $I$.

\section{Canonical module}
If $R$ is Cohen-Macaulay (which is known for the cases $n=2,3,4$) then its canonical module
is defined as
 $$\omega_{R}:={\rm Ext}_S^d(S/I,S)$$ 
 where $d=n^2-n$ is the height of $I$. Let $J=j_1,\ldots, j_{n^2-n}$ be the subideal of $I$ consisting of
 the off-diagonal elements in $XY-YX$. Then we have
\[
{\rm Ext}_S^d(S/I,S) \cong {\rm Ext}_{S/j_1}^d(S/I,S/j_1) \cong \cdots 
\cdots \cong {\rm Hom}_{S/J}(S/I,S/J) \cong (J:I)/J
\]
In the previous section we gave a conjecture on the generators of $(J:I)$. For $n=4$ we compute
$(J:I)$ using the conjecture and partially resolve $(J:I)/J$ using {\sc Macaulay}. 
We get the 
following Betti numbers:

\begin{verbatim}
% 1% 2% betti cp
total:     14   200   660  3821  
--------------------------------
    4:      3     -     -     -   
    5:      4   110   256    90   
    6:      7    90   908  3656 
    7:      -     -     6    75
\end{verbatim}

\noindent
This gives 
us the Betti numbers of the tail of the resolution of $I$ (see e.g. cor. 3.3.9 in~\cite{bh}). So we can compare this with
the Hilbert series of $S/I$:
\begin{eqnarray*}
h_{S/I}(t) &=& (1-15 t^2
     + 2 t^3
    +108 t^4
    -26 t^5
   -562 t^6
    +466 t^7
   +1613 t^8
  -2742 t^9
  -1078 t^{10}
   +5994 t^{11} \\
&&  -4367 t^{12}
  -2262 t^{13}
   +5630 t^{14}
  -3650 t^{15}
    +818 t^{16}
    +166 t^{17}
   -103 t^{18}
      +4 t^{19}
      +3 t^{20})
/(1-t)^{32}
\end{eqnarray*}
We see that our conjecture fits with the (last 6) coefficients of the polynomial in the numerator.
\smallskip
\noindent
Partly computing the resolution of $I$ we get the Betti numbers:
\begin{verbatim}

o18 = total: 1 16 115 595 2127 2791 848 60 5
          0: 1  .   .   .    .    .   .  . .
          1: . 15   2   .    .    .   .  . .
          2: .  . 108  30    3    .   .  . .
          3: .  .   4 565  466   45   4  . .
          4: .  .   .   . 1658 2746 844 60 5
\end{verbatim}
Splicing together these 2 Betti tables and using the Hilbert series we get
we get the following conjecture on the Betti numbers:

\[
\begin{array}{rrrrrrrrrrrrrr}
 total:& 1 &15& 115 &595 &2127& 4713& 6902+& 4432+& 5710+& 3821& 1170 &200 & 14 \\
          0: & {\bf 1} & . &  . &  . &   . &   . &   . &  . & .&.&.&.&. \\ 
          1: &. &{\bf 15} &  {\bf 2} &  . &   . &   . &   . & . &. &.&.&.&.\\
          2: &. & .& {\bf 108} & {\bf 30} &  {\bf 3} &   . &   . & . &. &.&.&.&.\\
          3: &. & .&   {\bf 4}& {\bf 565} & {\bf 466} &  {\bf 45} &   {\bf 4} & . &. &.&.&.&.\\
          4: &. & .&   .&   .& {\bf 1658}& {\bf 2746}&    {\bf 844}  &    {\bf 60} & {\bf 5} & .&.&.&.\\ 
          5: &. & .&   .&   .&    .&    1922& 6054 &  4372 &   c&   {\bf 75} &   {\bf 6}&.&.\\
          6: &. & .&   .&   .&    .&    .&    . &      d& 5705& {\bf 3656}& {\bf 908}&  {\bf 90}& {\bf 7} \\
          7: &. & . &  .&   .&    .&    .&    . &     . &   . &  {\bf 90}& {\bf 256} &{\bf 110} &{\bf 4} \\
          8: &. & . &  .&   .&    .&    .&    . &     . &   . &   .&   .&   .& {\bf 3} 
\end{array}
\]
where $-d+c = -2262$ (from the Hilbert series). The boldfaced numbers are the 2 earlier
macaulay computations and the others are based on the Hilbert series.

\section{Resolution}
Computer calculations indicate that there might be a (non-trivial) multiplicative structure 
on the resolution. The first $n-1$ lines in the Betti table can be interpreted as products
of the generators and the first syzygies. There is also an interesting "multiplicative pattern"
on the top "staircase" of the Betti numbers, 
i.e. these seem to be products 
of the $2$ linear first syzygies and the generators of the ideal. Consider below the Betti numbers for $n=3$
and partial the Betti numbers for $n=4,5,6$:

\newpage
\begin{verbatim}
n=3                                   n=4
o9 = total: 1 9 34 60 61 32 5         o18 = total: 1 16 115 595 2127 2791 848 60 5      
         0: 1 .  1  .  .  . .                   0: 1  .   1   .    .    .   .  . .  
         1: . 9  2  .  .  . .                   1: . 16   2   .    .    .   .  . .
         2: . . 31 32  3  . .                   2: .  . 108  30    3    .   .  . .
         3: . .  . 28 58 32 4                   3: .  .   4 565  466   45   4  . .
         4: . .  .  .  .  . 1                   4: .  .   .   . 1658 2746 844 60 5

\end{verbatim}
\begin{verbatim}
n=5                                         n=6
o9 = total: 1 25 291 2486 561 72            o14 = total: 1 36 605 6720 1199 105 4          
         0: 1  .   1    .   .  .                      0: 1  .   1    .    .   . .       
         1: . 25   2    .   .  .                      1: . 36   2    .    .   . .
         2: .  . 279   48   3  .                      2: .  . 598   70    3   . .
         3: .  .   4 2096 558 72                      3: .  .   4 6650 1196 105 4
         4: .  .   5  342

\end{verbatim}
\noindent
We seem to get $2(n^2-1)$ quadratic second syzygies, $3(n^2-1)$ quadratic fourth syzygies etc. We get
$2$ linear first syzygies, $3$ linear third syzygies, $4$ linear fifth syzygies etc. In the fourth line
of the table for $n=5$ we have $2096=\binom{24}{3} + 3 \cdot 24$ and $558=2 \cdot 279$.
We interpret the first 4 lines of this table as products,
let $f_1, \ldots, f_{24}$ be generators of the ideal,
                $g_1$, $g_2$ the linear first syzygies,
$h_1,h_2,h_3$ the first syzygies (that are not the Koszul relations) of degree $2$
and
$k_1, k_2, k_3, k_4$ the first syzygies of degree $3$, then the partial 
Betti table may be interpreted as follows:

\begin{small}
\begin{center}
\begin{tabular}{l|l|l|l|l|l}
hd $1$   & hd $2$  & hd $3$  &  hd $4$ & hd $5$ & hd $6$  \\ \hline
$|\{f_i\}|=24$ & $g_1, g_2$ & - & - & - & - \\ \hline
- & $|\{f_i  f_j\}| = \binom{24}{2}=276$& $|\{g_i  f_j\}|=48$  & $|\{g_i g_j\}|=3$ & - & - \\
 &  $h_1, h_2, h_3 $ & & & &  \\ \hline
& $k_1, k_2, k_3, k_4$ & $|\{f_i f_j f_k\} | = \binom{24}{3}$ & $|\{g_k f_i f_j\}| = 552$ & $|\{f_i g_j g_k\}| = 72$  & $|\{g_i g_j g_k\} | = 4$ \\
&  & $=2024$ & & & \\
&          &  $|\{f_i h_j\}| = 72 $      & $ |\{g_i h_j\}| = 6$ & & \\ \hline
\end{tabular}
\end{center}
\end{small}

\noindent
So up to a certain row (probably row $n-1$)
the generators and the first syzygies seem to generate everything (and the "multiplication" is nonzero). 
Our conjecture is that we have the
following Betti numbers for a general $n$:
\begin{small}
\begin{center}
\begin{tabular}{rllllllllccclllll}
\!\!\!\!\!\!\!\!\!\!\!\!&hd $1$   &
hd $2$  & hd $3$  &  hd $4$ & hd $5$ & hd $6$& hd $7$&
hd $8$ & $\cdots$& \!\!\!\!$n^2-n-1$ & \!\!\!\!$n^2-n$ \\ \hline
\!\!\!\!\!\!\!\!\!\!\!\!1: &$ n^2-1$ & $ 2$ & - & - & - & - & - & - & $\cdots$ &\!\!\!\!- &\!\!\!\!- \\ 
\!\!\!\!\!\!\!\!\!\!\!\!2:  &- & $ \binom{n^2-1}{2}+3$& $ 2(n^2-1)$ & $ 3$  & -
& - & - & -& $\cdots$ &\!\!\!\!-&\!\!\!\!- &\\ 
\!\!\!\!\!\!\!\!\!\!\!\!3:&-& $4$ &  p & p & $3 \cdot(n^2-1)$  & $4$ & - & -& $\cdots$ &\!\!\!\!- & \!\!\!\!- &\\
\!\!\!\!\!\!\!\!\!\!\!\!4:&-& $5$ &p&p&p&p & $\ddots$ & $\ddots$ &$\cdots$ &\!\!\!\!- &\!\!\!\!-\\
\!\!\!\!\!\!\!\!\!\!\!\!\vdots &-&-& \vdots &\vdots &\vdots&\vdots&\vdots&\vdots&$\ddots$ &\!\!\!\!-& \!\!\!\!-\\ 
\!\!\!\!\!\!\!\!\!\!\!\!$n-2$ &-& $n-1$ &p&p&p&p &p&p & $\cdots$ &\!\!\!\!-&\!\!\!\!-\\ 
\!\!\!\!\!\!\!\!\!\!\!\!$n-1$ &-& $n$ &?&?&?&?&?&? & $\cdots$ &\!\!\!\!-&\!\!\!\!-\\ 
\!\!\!\!\!\!\!\!\!\!\!\!$\vdots$ &-&-& \vdots &\vdots &\vdots&\vdots&\vdots&\vdots&\vdots & \!\!\!\!-&\!\!\!\!-\\ 
\!\!\!\!\!\!\!\!\!\!\!\!\!\!\!\!$\frac{n(n-1)}{2}$ &-&-& \vdots &\vdots &\vdots&\vdots&\vdots&\vdots&$\cdots$
&\!\!\!\!$\frac{n(n-1)}{2}(n^2-1)$ &\!\!\!\!$\frac{n(n-1)}{2} + 1$ \\
\!\!\!\!\!\!\!\!\!\!\!\!\vdots &-&-& \vdots &\vdots &\vdots&\vdots&\vdots&\vdots&\vdots&\!\!\!\!\vdots
&\!\!\!\!\vdots\\ 
\!\!\!\!\!\!\!\!\!\!\!\!$M$:&-&-&$\cdot$&$\cdot$&$\cdot$&$\cdot$&$\cdot$&$\cdot$&
$\cdot$&$\!\!\!\!\cdot$&\!\!\!\!\!\!\!\!$s(k-s+1)+1$ \\
\!\!\!\!\!\!\!\!\!\!\!\!$M+1$:&-&-&-&-&-&-&-&-&-&\!\!\!\!-&\!\!\!\!- \\
\end{tabular}
\end{center}
\end{small}

\smallskip
\noindent
where $p$ means products of earlier entries and $M$ is determined by $d_{min}$ from 
the conjecture on the canonical module i.e. $M=n(n-1) - d_{min}$. The numbers $s$ and $k$
are defined in the section on $(J:I)$.

\section{Minimal primes of $J$}

It has been conjectured that $I$ is a prime ideal and this can be deduced from
$I$ being Cohen-Macaulay~\cite{vasc} which we know for $n=2$, $n=3$ and $n=4$ so at least in those cases
$I$ is a minimal prime of $J$. Knutson~\cite{kn} gives a conjecture that $J$ has one other minimal
prime which is generated by determinants of matrices whose columns are the diagonals of powers of $X$
and $Y$. These determinants are elements of $(J:I)$ but, if our conjecture is true, do not generate it for all $n$. Our conjecture is that the other prime ideal
is given by $(J:I)$ which coincides with Knutson's equations for $n=3$. For $n=4$ we get by picking 
the syzygies $(E,X,Y,XY+YX)$ an element of bidegree $(2,2)$ that is a zero-divisor on $J$ and
that can not, for bidegree reasons, be created by determinants coming from the diagonals of the powers of $X$ and $Y$.

\bigskip
\noindent
Science Institute, University of Iceland, Dunhagi 3, IS-107 Reykjavik, Iceland \\
{\em E-mail address:} freyjah@raunvis.hi.is


\begin{thebibliography}{99}


 \bibitem{bast} D. Bayer, M. Stillman, {\em Macaulay: A system for
computation in algebraic geometry and commutative algebra.} Source and object
 code available for Unix and Macintosh computers. Contact the authors or
 download from zariski.harvard.edu via anonymous ftp. (1990)

 \bibitem{bss} D.
 Bayer, M. Stillman, Ma. Stillman,         {\ Macaulay User Manual.}
 
 
\bibitem{bh} W. Bruns, J. Herzog, {\em Cohen-Macaulay rings,} Cambridge University Press, 1993. 


\bibitem{for} E. Formanek, {\em The Polynomial Identities and Invariants
       of $n \times n$ matrices,}  CBMS Regional Conference              Series
 in Mathematics {\bf 78}, published for the Conference              Board of
 the Mathematical             Sciences, Washington, DC, (1991).

\bibitem{gerst}
 M. Gerstenhaber, {\em On dominance and varieties of                 commuting
 matrices,}  Ann. of Math. {\bf 73} (1961), 324-348.

\bibitem{m2} D. Grayson, M. Stillman, Macaulay 2: a computer algebra system 
for algebraic geometry and commutative algebra, available at http://www.math.uiuc.edu/Macaulay2.
 
 \bibitem{frh1}
 F. Hreinsdottir, {\em A Case Where Choosing a Product Order Makes the
   Calculations of a Groebner Basis Much Faster,}           J. Symbolic
 Comput. {\bf 18} (1994), 373-378.


\bibitem{frh4} F. Hreinsdottir, {\em On the ring
 of Commuting Matrices,} thesis Stockholm University 1997.


 \bibitem{frh2} F. Hreinsdottir, {\em The Koszul
 Dual of the Ring of             Commuting Matrices,} Comm. Algebra {\bf 26} (1998), 3807-3819.

\bibitem{kn} A. Knutson, {\em Some Schemes Related to the Commuting Variety}, to appear in J. Algebraic Geom.,
ArXiv: math.AG/0306275, 2003. 

 \bibitem{mot} T. Motzkin
 and O. Taussky, {\em Pairs of matrices with                property L II,}
 Trans. Amer. Math. Soc. {\bf 80} (1955),                 387-401.


 \bibitem{vasc} W. V. Vasconcelos, {\em Computational
 Methods in Commutative Algebra and Algebraic Geometry,} Algorithms and Computation in Math. {\bf 2}, 
 Springer 1998.
 
\end{thebibliography}
\end{document}